\documentclass{amsart}
\usepackage{pgf}
\usepackage{amssymb}
\usepackage{dsfont}
\usepackage{txfonts}
\usepackage{graphicx}
\usepackage{amsmath}
\usepackage[colorlinks=black]{hyperref}

\newtheorem{theo}{Theorem}[section]
\newtheorem{definition}[theo]{Definition}
\newenvironment{defi}{\begin{definition}\rm}{\end{definition}}
\newtheorem{remarque}[theo]{Remark}
\newenvironment{remark}{\begin{remarque}\rm}{\end{remarque}}
\newtheorem{exemple}[theo]{Example}
\newenvironment{ex}{\begin{exemple}\rm}{\end{exemple}}
\newtheorem{lemma}[theo]{Lemma}

\newtheorem{coro}[theo]{Corollary}
\newtheorem{nota}[theo]{Notation}

\newtheorem{prf1}{\it {Idea of proof}}

\newtheorem{prf}{\it {Proof}}

\newenvironment{demo}{\begin{prf}\rm}{\hfill$\Box$\end{prf}}

\def\varLim@#1#2{%
\vtop{\m@th\ialign{##\cr
\hfil$#1\operator@font Lim$\hfil\cr
\noalign{\nointerlineskip\kern1.5\ex@}#2\cr
\noalign{\nointerlineskip\kern-\ex@}\cr}}%
}
\def\varinjLim{%
\mathop{\mathpalette\varLim@{\rightarrowfill@\textstyle}}\nmlimits@
}

\title{On the algebraicity of Puiseux series.} 

\author{Michel Hickel and Micka\"{e}l Matusinski}

\address{Universit\'{e} de Bordeaux\\
IMB Institut de Math\'{e}matiques de Bordeaux\\
351 cours de la Lib\'{e}ration\\
33405 Talence Cedex\\
France}

\email{michel.hickel@math.u-bordeaux1.fr\\
mickael.matusinski@math.u-bordeaux1.fr}

\begin{document}

\begin{abstract}
We deal with the algebraicity of a Puiseux series in terms of the properties of its coefficients. We show that the algebraicity of a Puiseux series for given bounded degree is determined by a finite number of explicit polynomial formulae. Conversely, given a vanishing polynomial, there is a closed-form formula for the coefficients of the series in terms of the coefficients of the polynomial and of an initial part of the series. 
\end{abstract}
\maketitle

\section{Introduction.}

Let  $K$ be a zero characteristic field and $\overline{K}$ its algebraic closure. We consider  $K[[x]]$, the domain of formal power series with coefficients in $K$, and its fraction field $K((x))$. We denote by $K((\hat{x})):= \displaystyle\bigcup_{n=1}^{\infty}K((x^{1/n}))$ the field of  formal Puiseux series (with coefficients in $K$).
 By the Newton-Puiseux theorem (see e.g. \cite[Theorem 3.1]{walker_alg-curves} and \cite[Proposition p.314 ]{rib-vdd_ratio-funct-field}),  an algebraic closure of $K((x))$ is given by  $\mathcal{P}_K:=\displaystyle\bigcup_L L((\hat{x}))$ where $L$ ranges over the finite extensions of $K$ in $\overline{K}$. In particular, if  $K=\overline{K}$, then $\mathcal{P}_K=K((\hat{x}))$.   Among Puiseux series, we are interested in algebraic ones, say the Puiseux series which verify a polynomial equation with coefficients that are themselves polynomials in $x$:  $P(x,y)\in K[x][y]$.\\

Among the numerous works concerning algebraic Puiseux series \cite{vdpoorten:alg-series,flajolet-sedgewick:ana-combinatorics,banderier-drmota:coeff-alg-series}, we deal with the following questions:\\
\noindent$\bullet$ \textbf{Reconstruction of a vanishing polynomial for a given algebraic Puiseux series.} Generically, a vanishing polynomial of a given algebraic power series can be computed as a Pad\'{e}-Hermite approximant \cite[Chap. 7]{salvy-et-al:cours-mpri}. In fact, the algebraicity of a Puiseux series can be encoded by the vanishing of certain determinants derived from the coefficients of the series. We extend this approach by showing how to reconstruct the coefficients of a vanishing polynomial by means of some minors of these determinants (see Section \ref{section:Wilc}). More precisely, we show that, for a given bounded degree, there are finitely many universal polynomials allowing to check the algebraicity of a series and to perform this reconstruction (see Theorem \ref{theo:wilc}). \\
\noindent$\bullet$ \textbf{Description of the  coefficients of an algebraic Puiseux series in terms of the coefficients of a vanishing polynomial.} An approach consists in considering that the series coefficients verify a linear recurrence relation, which allows an asymptotic computation of the coefficients. This property follows classically from the fact that an algebraic Puiseux series is \textit{differentiably finite (D-finite)}, that is, it satisfies a linear differential equation with coefficients in $K[x]$ \cite{comtet:coeff-alg-series,stanley:generating-funct,stanley:enum-combinatorics,singer:lin-ODE-alg-funct,chudnovsky-s:alg-series-I,chudnovsky-s:alg-series-II,bostan-et-al:diff-equ-alg-funct}. 

Another approach consists in determining a closed-form expression in terms of the coefficients of a vanishing polynomial. In this direction, P. Flajolet and M. Soria (see the habilitation thesis of M. Soria (1990) and \cite{flajolet-soria:coeff-alg-series}) proposed a formula in the case of a series satisfying a reduced Henselian equation (see the Definition \ref{defi:equ-hensel-red} for this terminology) with complex coefficients. This formula extends to coefficients in an arbitrary zero characteristic field $K$ via a work of Henrici \cite{henrici:lagr-burmann}.

Here we complete this approach to the case of a Puiseux series which satisfies a general polynomial equation $P(x,y)=0$, by showing that the coefficients of such series can be computed applying the Flajolet-Soria formula to a  polynomial  naturally derived from $P$ (see Section \ref{section:flajo-soria}).  \\

\noindent\textbf{Acknowledgement}. The authors would like to thank kindly G. Rond for his observation that allowed us to simplify the proof of the Theorem \ref{theo:wilc}. We are also indebted to K. Kurdyka and O. Le Gal for valuable discussions. 

\section{Preliminaries} 
Let us denote $\mathbb{N}:=\mathbb{Z}_{\geq 0}$ and $\mathbb{N}^*:=\mathbb{N}\setminus\{0\}=\mathbb{Z}_{>0}$. For any set $\mathcal{E}$, we will write $|\mathcal{E}|:=\mathrm{Card}(\mathcal{E})$. For any vector of natural numbers $K=(k_1,\ldots,k_n)$, we set $K!:=\displaystyle\prod_{i=1}^nk_i!$, $|K|:=\displaystyle\sum_{i=1}^nk_i$ and $  \|K\| :=\displaystyle\sum_{i=1}^ni\,{k_i}$. The floor function will be written $\lfloor q \rfloor$ for $q\in\mathbb{Q}$. \\

Let  $\tilde{y}_0=\displaystyle\sum_{n\geq n_0} \tilde{c}_nx^{n/p}\in K((\hat{x})),\ \tilde{c}_{n_0}\neq 0$, a Puiseux series. We denote $$\tilde{y}_0=x^{(n_0-1)/p}\displaystyle\sum_{n\geq 1} c_nx^{n/p}=x^{(n_0-1)/p}\tilde{z}_0\ \ \textrm{with}\  c_1\neq 0.$$ 
The series $\tilde{y}_0$ is a root of a polynomial $\tilde{P}(x,y)=\displaystyle\sum_{i,j} \tilde{a}_{i,j}x^iy^j$ of degree $d_y$ in $y$ if and only if the series $y_0=\displaystyle\sum_{n\geq 1} c_nx^n$ is a root of  $x^m\tilde{P}(x^p,x^{n_0-1}y)$, the latter being a polynomial for $m=\max\{0;(1-n_0)d_y\}$. The existence of a polynomial $\tilde{P}$ cancelling $\tilde{y}_0$ is equivalent to the one of a polynomial $P(x,y)=\displaystyle\sum_{i,j} a_{i,j}x^iy^j$ cancelling $y_0$, such that, for $(i,j)$ belonging to the support of $P$, one has $i\equiv (n_0-1)j\ [p]$ if $n_0-1\geq 0$, and $i\equiv (1-n_0)(d_y-j)\ [p]$ otherwise. Thus the algebraicity of $\tilde{y}_0$ is equivalent to that of $y_0$ but \textit{ with constraints}. This leads us to the following definition:

\begin{defi}\label{defi:alg-relative}
Let $\mathcal{F}$ and $\mathcal{G}$ be two strictly increasing finite sequences of  couples $(i,j)\in\mathbb{N}^2$ ordered anti-lexicographically: 
$$(i_1,j_1) \leq (i_2,j_2)\Leftrightarrow  j_1 < j_2\textrm{ or } (j_1 = j_2\ \textrm{et}\ i_1 \leq i_2).$$
We suppose additionally that  $\mathcal{F}\geq (0,1)>\mathcal{G}>(0,0)$ (thus the elements of $\mathcal{G}$ are couples of the form $(i,0)$, $i\in\mathbb{N}^*$, and those of  $\mathcal{F}$ are of the form  $(i,j),\ j\geq 1$). We say that a series $y_0=\displaystyle\sum_{n\geq 1} c_nx^n\in K((x))$, $c_1\neq 0$, is \textbf{algebraic relatively to $(\mathcal{F},\mathcal{G})$} if there exists a polynomial $P(x,y)=\displaystyle\sum_{(i,j)\in\mathcal{F}\cup\mathcal{G}} a_{i,j}x^iy^j\in K[x,y]$ such that $P(x,y_0)=0$.
\end{defi}

Flajolet and Soria (see the habilitation thesis of M. Soria (1990) and \cite{flajolet-soria:coeff-alg-series}) gave a closed-form expression to compute the coefficients of a formal solution of a reduced Henselian equation in the following sense:

\begin{defi}\label{defi:equ-hensel-red}
We call \textbf{reduced Henselian equation} any equation of the following type:
$$y=Q(x,y)\ \textrm{ with }\ Q(x,y)\in K[x,y],$$
such that $Q(0,0)=\displaystyle\frac{\partial Q}{\partial y}(0,0)=0$ and $Q(x,0)\nequiv 0$.
\end{defi}

\begin{theo}[Flajolet-Soria formula]\label{theo:formule-FS}
Let $y=Q(x,y)=\displaystyle\sum_{i,j}a_{i,j}x^iy^j$ be a reduced Henselian equation. Then the  coefficients of the unique solution $\displaystyle\sum_{n\geq 1}c_nx^n$ are given by:
$$ c_n=\displaystyle\sum_{m=1}^{2n-1}\frac{1}{m}\displaystyle\sum_{|\underline{k}|=m,\ ||\underline{k}||_1=n,\ ||\underline{k}||_2=m-1}\frac{m!}{\prod_{i,j}k_{i,j}!}\prod_{i,j}a_{i,j}^{k_{i,j}},$$
where $\underline{k}=\displaystyle(k_{i,j})_{i,j}$, $\ |\underline{k}|=\displaystyle\sum_{i,j}k_{i,j}$, $\ ||\underline{k}||_1 = \displaystyle\sum_{i,j}i\, k_{i,j}$ and $\ ||\underline{k}||_2 = \displaystyle\sum_{i,j}j\, k_{i,j}$.
\end{theo}

\begin{remark}\label{rem:FS}
Let us consider the particular case where the coefficients of $Q$ verify $a_{0,j}=0$ for all $j$.  So, for any $\underline{k}$ such that $|\underline{k}|=m$ and $\displaystyle\prod_{i,j}a_{i,j}^{k_{i,j}}\neq 0$, we have that $\|\underline{k}\|_1\geq m$. Thus, to have $\|\underline{k}\|_1=n$, one needs to have $m\leq n$. The Flajolet-Soria formula can be written:
$$ c_n=\displaystyle\sum_{m=1}^{n}\frac{1}{m}\displaystyle\sum_{|\underline{k}|=m,\ ||\underline{k}||_1=n,\ ||\underline{k}||_2=m-1}\frac{m!}{\prod_{i,j}k_{i,j}!}\prod_{i,j}a_{i,j}^{k_{i,j}}.$$
\end{remark}

\section{Characterizing the algebraicity of a formal power series}\label{section:Wilc}
Here we resume the results from  \cite{wilczynski:alg-power-series}.
Suppose we are given a series $y_0=\displaystyle\sum_{n\geq 1} c_nx^{n}\in K((x))$ with  $c_1\neq 0$. For any $j\in\mathbb{N}$, consider the multinomial expansion of $y_0^j$, that we denote:
$$y_0^j=\displaystyle\sum_{n\geq 1} c_n^{(j)}x^{n}.$$
Of course, one has that $c_n^{(j)}=0$ for $n<j$ and $c_j^{(j)}=c_1^j\neq 0$. For $j=0$, let $y_0^0=1$. We remark that for any $n$ and any $j$, $c_n^{(j)}$ is a homogeneous polynomial with natural number coefficients of degree $j$ in the $c_m$ for $m\leq n-j+1$. 
\begin{defi}\label{defi:mat_Wilc}
\begin{enumerate}
    \item 
Given a couple $(i,j)\in\mathbb{N}\times\mathbb{N}$, we call \textbf{Wilczynski vector} $V_{i,j}$ the infinite vector with components:\\
- if $j\geq 1$, a sequence of $i$ zeros followed by the coefficients of the multinomial expansion $y_0^j$:
$$V_{i,j}:=(0,\ldots,0,c_1^{(j)},c_2^{(j)},\ldots,c_n^{(j)},\ldots);$$
- otherwise, 1  in the $i$th position and 0 for the other coefficients
$$V_{i,0}:=(0,\ldots,1,0,0,\ldots,0,\ldots).$$
\item Let $\mathcal{F}$ and $\mathcal{G}$ be two sequences as in the Definition \ref{defi:alg-relative}. We associate to $\mathcal{F}$ and $\mathcal{G}$ the \textbf{ (infinite)  Wilczynski matrix } whose columns are the corresponding vectors $V_{i,j}$:
$$M_{\mathcal{F},\mathcal{G}}:=(V_{i,j})_{(i,j)\in\mathcal{F}\cup \mathcal{G}}\,.$$
We define also the \textbf{reduced Wilczynski matrix}, $M_{\mathcal{F},\mathcal{G}}^{red}$: it is the matrix obtained from $M_{\mathcal{F},\mathcal{G}}$ by removing the columns indexed in $\mathcal{G}$,  and also removing the corresponding rows (suppress the $i$th row for any $(i,0)\in\mathcal{G}$). This amounts exactly to remove the rows containing the coefficient 1 for some Wilczynski vector indexed in $\mathcal{G}$. 
\end{enumerate}
\end{defi}

\begin{lemma}[Wilczynski]\label{lemme:wilcz}
The series $y_0$ is algebraic relatively to $(\mathcal{F},\mathcal{G})$ if and only if all the minors of order $|\mathcal{F}\cup\mathcal{G}|$ of the  Wilczynski matrix $M_{\mathcal{F},\mathcal{G}}$ vanish, or also if and only if all the minors of order  $|\mathcal{F}|$ of the reduced Wilczynski matrix $M_{\mathcal{F},\mathcal{G}}^{red}$ vanish.
\end{lemma}
\begin{demo}
Given a   non trivial polynomial $P(x,y)=\displaystyle\sum_{(i,j) \in\mathcal{F}\cup\mathcal{G}}a_{i,j}x^iy^j$,  we compute:
$$\begin{array}{lcl}
 P(x,y_0)&=&\displaystyle\sum_{(i,j) \in\mathcal{F}}a_{i,j}x^i\left(\displaystyle\sum_{n\geq 1} c_n^{(j)}x^{n}\right)+\displaystyle\sum_{(i,0) \in\mathcal{G}}a_{i,0}x^i\\
&=& M_{\mathcal{F},\mathcal{G}}\cdot (a_{i,j})_{(i,j)\in\mathcal{F}\cup\mathcal{G}}\,.
\end{array}$$
where the components of the infinite vector thus obtained are the coefficients of the expansion of  $P(x,y_0)$ with respect to the powers of $x$ in increasing order. The series $y_0$ is a root of $P$ if and only if this infinite vector is the zero vector, which means that the rank of  $M_{\mathcal{F},\mathcal{G}}$ is less than  $|\mathcal{F}\cup\mathcal{G}|$, the number of columns of  $M_{\mathcal{F},\mathcal{G}}$. The latter condition is characterized as in finite dimension by the vanishing of all the minors of maximal order.\\ 
Let us now remark that, in the infinite vector $M_{\mathcal{F},\mathcal{G}}\cdot (a_{i,j})_{(i,j)\in\mathcal{F}\cup\mathcal{G}}$, if we remove the components of number $i$ for $(i,0)\in\mathcal{G}$, then we get exactly the infinite vector $M_{\mathcal{F},\mathcal{G}}^{red}\cdot (a_{i,j})_{(i,j)\in\mathcal{F}}$. The vanishing of the latter means precisely that the rank of  $M_{\mathcal{F},\mathcal{G}}^{red}$ is less than $|\mathcal{F}|$. 
Conversaly, if the columns of  $M_{\mathcal{F},\mathcal{G}}^{red}$ are dependent for certain $\mathcal{F}$ and $\mathcal{G}$, we denote by $(a_{i,j})_{(i,j)\in\mathcal{F}}$ a corresponding sequence of coefficients of a non trivial vanishing linear combination of the column vectors. Then it suffices to note that the remaining  coefficients $a_{k,0}$ for  $(k,0)\in\mathcal{G}$ are each uniquely determined  as follows: 
\begin{equation}\label{equ:terme-cst}
a_{k,0}=-\displaystyle\sum_{(i,j)\in\mathcal{F}, i< k} a_{i,j}c_{k-i}^{(j)}\,.
\end{equation} 
\end{demo}

We deal with the implicitization problem for algebraic power series: for fixed bounded degrees in $x$ and $y$,  given the expression of an algebraic series, can we reconstruct a vanishing polynomial ? if yes, how ?

\begin{defi}\label{defi:poly-wilc}
Let us consider a reduced Wilczynski matrix $M_{\mathcal{F},\mathcal{G}}^{red}$ associated to two sequences $\mathcal{F}$ and $\mathcal{G}$ of couples $(i,j)$ as in \ref{defi:alg-relative}. We call \textbf{Wilczynski polynomial} any  polynomial in the  coefficients $c_n$ of $y_0$ obtained as a minor of $M_{\mathcal{F},\mathcal{G}}^{red}$. We denote such  Wilczynski polynomial by $Q_{\underline{k},\underline{I}}$, where  $\underline{I}:=((i_1,j_1),\ldots,(i_l,j_l))$ is a subsequence of  $\mathcal{F}$ indicating the $l$ columns of $M_{\mathcal{F},\mathcal{G}}^{red}$, and  $\underline{k}:=(k_1,k_2,\cdots,k_l)$ a strictly increasing sequence of natural numbers indicating the $l$ rows of  $M_{\mathcal{F},\mathcal{G}}^{red}$ used to form the minor of $M_{\mathcal{F},\mathcal{G}}^{red}$. One  has that $l\in\mathbb{N}^*$, $l\leq |\mathcal{F}|$, $l$ being the order of that minor, that we will also call the  \textbf{order} of the Wilczynski polynomial $Q_{\underline{k},\underline{I}}$. Note also  that a  Wilczynski polynomial $Q_{\underline{k},\underline{I}}$ is  homogeneous of degree equal to $\displaystyle\sum_{(i,j)\in \underline{I},\ c_k^{(j)}\nequiv 0}j\ $ (indeed, the coefficients of $M_{\mathcal{F},\mathcal{G}}^{red}$ verify: $c_k^{(j)}\equiv 0\Leftrightarrow k<j$). By convention, we call  \textbf{Wilczynski polynomial of order 0} any non zero constant polynomial. 
\end{defi}

By  \ref{lemme:wilcz}, the algebraicity of $y_0$ for certain $\mathcal{F}$ and $\mathcal{G}$ is equivalent to the vanishing of all the  $Q_{\underline{k},\mathcal{F}}$ of order $l=|\mathcal{F}|$, for the specific values of the  given $c_n$, coefficients of $y_0$. \\

\begin{ex}\label{ex:wilc} Let $y_0=\displaystyle\sum_{n\geq 1} c_nx^{n}\in K((x))$ be a series with  $c_1\neq 0$.
We consider the conditions for $y_0$ to be a root of a polynomial of type:
$$P(x,y)=a_{2,0}x^2+a_{2,1}x^2y+(a_{0,2}+a_{2,2}x^2)y^2.$$
Thus, $\mathcal{F}=\{(2,1),(0,2),(2,2)\}$ and $\mathcal{G}=\{(2,0)\}$. The corresponding Wilczynski matrix is:
$$ M :=\left[\begin{array}{cccc}0 & 0 & 0 & 0 \\1
 & 0 & \mathop{\rm c}_{1}^{2} & 0 \\0 & \mathop{\rm c}_{1} & 2\cdot c 
_{1}\cdot c _{2} & 0 \\0 & \mathop{\rm c}_{2} & c _{2}^{2}+2~c _{1}~c 
_{3} & \mathop{\rm c}_{1}^{2} \\0 & \mathop{\rm c}_{3} & 2~c _{1}~c _{
4}+2~c _{2}~c _{3} & 2\cdot c _{1}\cdot c _{2} \\0 & \mathop{\rm c}_{4
} & 2~c _{2}~c _{4}+c _{3}^{2}+2~c _{1}~c _{5} & c _{2}^{2}+2~c _{1}~c
 _{3} \\ \vdots & \vdots & \vdots & \vdots \\\end{array}\right],$$
and the reduced matrix is:
$$M ^{{\it red} }:=
\left[\begin{array}{ccc}0 & 0 & 0 \\c _{1} & 2\cdot c _{1}\cdot c _{2}
 & 0 \\c _{2} & c _{2}^{2}+2~c _{1}~c _{3} & \mathop{\rm c}_{1}^{2} \\
c _{3} & 2~c _{1}~c _{4}+2~c _{2}~c _{3} & 2\cdot c _{1}\cdot c _{2}
 \\c _{4} & 2~c _{2}~c _{4}+c _{3}^{2}+2~c _{1}~c _{5} & c _{2}^{2}+2~
c _{1}~c _{3} \\ \vdots & \vdots & \vdots \\\end{array}\right].$$
We give the four first non trivial  Wilczynski polynomials of maximal order 3, equals to  minors 3x3 of $M^{red}$. So one has that $\underline{I}=\mathcal{F}$ as index for $Q_{\underline{k},\underline{I}}$ :
\begin{center}$\begin{array} {rl} Q_{\underline{k},\mathcal{F}}:=&-2\,{c_{{1}}}^{2} \left( {c_{{2}}}^{3}-2\,c_{{3}}c_{{1}}c_{{2}}+{c_{{1}}}^{2}c_{{4}} \right)\ \textrm{ pour }\underline{k}=(2,3,4),  \\
Q_{\underline{k},\mathcal{F}}:=&-c_{{1}} \left( {c_{{2}}}^{4}-3\,{c_{{1}}}^{2}{c_{{3}}}^{2}+2\,{c_{{1}}}^{3}c_{{5}} \right)\ \textrm{ pour }\underline{k}=(2,3,5),   \\
Q_{\underline{k},\mathcal{F}}:=&-2\,{c_{{1}}}^{2} \left( -c_{{4}}{c_{{2}}}^{2}-2\,c_{{1}}c_{{4}}c_{{3}}+c_{{2}}{c_{{3}}}^{2}+2\,c_{{1}}c_{{2}}c_{{5}} \right)\ \textrm{ pour }\underline{k}=(2,4,5),   \\
Q_{\underline{k},\mathcal{F}}:=& 8\,c_{{2}}{c_{{1}}}^{2}c_{{4}}c_{{3}}+{c_{{2}}}^{4}c_{{3}}-2\,{c_{{2}}}^{2}{c_{{3}}}^{2}c_{{1}}-4\,{c_{{1}}}^{2}{c_{{2}}}^{2}c_{{5}}-3\,{c_{{1}}}^{2}{c_{{3}}}^{3}+2\,c_{{3}}{c_{{1}}}^{3}c_{{5}}-2\,{c_{{1}}}^{3}{c_{{4}}}^{2}\\
& \textrm{ pour }\underline{k}=(3,4,5)  .
\end{array}$\end{center}
The series $y_0$ is a root of a polynomial $P(x,y)$ as above if and only if all the Wilczynski polynomials of order 3 vanish. This implies in particular that:
$$c_{{4}}=-{\frac {c_{{2}} \left( {c_{{2}}}^{2}-2\,c_{{1}}c_{{3}}
 \right) }{{c_{{1}}}^{2}}}\ \textrm{ and }\ 
c_{{5}}=-\,{\frac {{c_{{2}}}^{4}-3\,{c_{{1}}}^{2}{c_{{3}}}^{2}}{{2\, c_
{{1}}}^{3}}}.$$

\end{ex}

\begin{theo}\label{theo:wilc}
Let $\mathcal{F}$ and $\mathcal{G}$ be two finite sequences of couples as in  \ref{defi:alg-relative}. We set  $d_y:=\max\{j,\ (i,j)\in\mathcal{F}\}$, $d_x:=\max\{i,\ (i,j)\in\mathcal{F}\cup\mathcal{G}\}$ and $N:=2d_xd_y$.  Then there exists a finite number of homogeneous polynomials  $a_{i,j}^{(\lambda)}\in\mathbb{Z}[x_1,\ldots,x_N]$, $(i,j)\in\mathcal{F}\cup\mathcal{G}$, $\lambda\in\Lambda$, of total degree $\deg a_{i,j}^{(\lambda)}\leq \frac{1}{2}d_y(d_y+1)(d_x+1)-1$ for $(i,j)\in\mathcal{F}$, and $\deg a_{i,0}^{(\lambda)}\leq \frac{1}{2}d_y(d_y+1)(d_x+1)-1+i$ for  $(i,0)\in\mathcal{G}$, such that, for any $y_0=\displaystyle\sum_{n\geq 1} c_nx^{n}\in K[[x]]$ with $c_1\neq 0$  series algebraic relatively to $(\mathcal{F},\mathcal{G})$, there is $\lambda\in\Lambda$ such that the polynomial:
$$ P^{(\lambda)}(x,y)=\displaystyle\sum_{(i,j) \in\mathcal{F}}a_{i,j}^{(\lambda)}(c_1,\ldots,c_N)x^iy^j+\displaystyle\sum_{(i,0) \in\mathcal{G}}a_{i,0}^{(\lambda)}(c_1,\ldots,c_N)x^i\ \in K[x,y]$$
vanishes at $y_0$. 
\end{theo}
First , we give the reconstruction process. Then we will show its finiteness.
\begin{demo}
Let $y_0=\displaystyle\sum_{n\geq 1} c_nx^{n}\in K[[x]]$ with $c_1\neq 0$ be algebraic relatively to $(\mathcal{F},\mathcal{G})$. We show how to reconstruct a vanishing polynomial   $P(x,y)$ of $y_0$.

Let $Q(x,y)=\displaystyle\sum_{(i,j) \in\mathcal{F}}b_{i,j}x^iy^j+\displaystyle\sum_{(i,0) \in\mathcal{G}}b_{i,0}x^i$ be a polynomial that vanishes at $y_0$. We proceed by induction on  $m$ the number of non zero coefficients $b_{i,j}$ for $(i,j)\in\mathcal{F}$. If $m=1$, $Q(x,y)$ is of the form:
$$Q(x,y)=b_{i,j}x^iy^j+\displaystyle\sum_{(i,0) \in\mathcal{G}}b_{i,0}x^i,$$
with $b_{i,j}\neq 0$. So we must have that $b_{n,0}=0$ for $n<i+j$, and the series $y_0$ verifies: 
$$\displaystyle\sum_{(n,0) \in\mathcal{G}}b_{n,0}x^n=-b_{i,j}x^iy_0^j=\displaystyle\sum_{n\geq i} -b_{i,j}c_{n-i}^{(j)}x^{n}.$$
The criterion \ref{lemme:wilcz} means here that the order 1 minors of  $M_{(i,j),\mathcal{G}}^{red}$,  being equal to $c_{n-i}^{(j)}$ for $(n,0)\notin \mathcal{G}$, are all null. We fix the coefficient $a_{i,j}$ arbitrarily in $\mathbb{Z}\setminus \{0\}$: it is a  constant Wilczynski polynomial. Then the other coefficients are uniquely determined in accordance with the relation (\ref{equ:terme-cst}) by the equation: 
$$a_{n,0}:=-a_{i,j}c_{n-i}^{(j)},\ \ (n,0)\in\mathcal{G}.$$
Thus the coefficient $a_{n,0}$ is a polynomial of degree $j$ in the $c_k$, $k\leq n-i-j+1$, which verifies indeed that $j\leq d_y \leq \frac{1}{2}d_y(d_y+1)(d_x+1)\leq \frac{1}{2}d_y(d_y+1)(d_x+1)-1+n$. \\
Consider now the case where the vanishing polynomial  $Q(x,y)$ of  $y_0$ has $m\geq 2$ non zero terms. So there exists a subfamily $\mathcal{F}'\subset \mathcal{F}$, the indices of the non zero coefficients of $Q(x,y)$, with $|\mathcal{F}'|=m$ and such that the order $m$ minors of $M_{\mathcal{F}',\mathcal{G}}^{red}$ all vanish. Suppose that there exists an order $m-1$ minor of this matrix, i.e. a  Wilczynski polynomial $Q_{\underline{k}_0,\underline{I}_0}$, which is non zero. Denote by $M_{\underline{k}_0,\underline{I}_0}$ the square matrix  whose determinant is this minor, and $C_{\underline{k}_0,(i_0,j_0)}$ the  $p_0$-th column that has been removed to form this minor. We get a Cramer system of equations: 
$$M_{\underline{k}_0,\underline{I}_0}\cdot  (b_{i,j})_{(i,j)\neq (i_0,j_0)}=-b_{i_0,j_0}C_{\underline{k}_0,(i_0,j_0)}.$$
Let us build the coefficients $a_{i,j}$ verifying:
$$M_{\underline{k}_0,\underline{I}_0}\cdot  (a_{i,j})_{(i,j)\neq (i_0,j_0)}=-a_{i_0,j_0}C_{\underline{k}_0,(i_0,j_0)},$$
by taking $a_{i_0,j_0}:=(-1)^{p_0}Q_{\underline{k}_0,\underline{I}_0}$ and by computing the other $a_{i,j}$ by Cramer's rule. Thus the $a_{i,j}$ are all order $m-1$ minors of $M_{\mathcal{F}',\mathcal{G}}^{red}$, and so, up to the sign, Wilczynski polynomials $ Q_{\underline{k}_0,\underline{I}}$ of order $m-1$.  If $\underline{k}_0=(k_{0,1},\ldots,k_{0,m-1})$, we set: 
\begin{equation}\label{equ:N} N_{y_0}:=k_{0,m-1}.\end{equation}
The $a_{i,j}$ are  homogeneous polynomials of $\mathbb{Z}[x_1,\ldots,x_{N_{y_0}}]$. The degree of a Wilczynski polynomial $Q_{\underline{k}_0,\underline{I}}$ verifies:
$$\begin{array}{lcl}
\deg Q_{\underline{k}_0,\underline{I}}&=&\displaystyle\sum_{(i,j)\in \underline{I},\ c_k^{(j)}\nequiv 0}j\\
&\leq& -1+\displaystyle\sum_{(i,j)\in\mathcal{F} }j\\
&\leq& -1+(d_x+1)\displaystyle\sum_{j=1}^{d_y}j\\
&=& \frac{1}{2}d_y(d_y+1)(d_x+1)-1.
\end{array}$$
The coefficients $a_{n,0}$ for $(n,0)\in\mathcal{G}$ are obtained via the relations (\ref{equ:terme-cst}):
$$a_{n,0}=-\displaystyle\sum_{(i,j)\in\mathcal{F}, n>i} a_{i,j}c_{n-i}^{(j)}.$$ 
Knowing that  $c_{n-i}^{(j)}\nequiv 0 \Rightarrow n-i\geq j$, and in this case $\deg c_{n-i}^{(j)}=j$, we deduce that $\deg a_{n,0}\leq n+\max_{(i,j)\in\mathcal{F}}(\deg a_{i,j})$ as desired. The polynomial $P(x,y)=\displaystyle\sum_{(i,j)\in\mathcal{F}'\cup\mathcal{G}}a_{i,j}x^iy^j$ is proportional to $Q$, so it vanishes at $y_0$.\\
Suppose now that all the minors of order $m-1$ are zero. So, restricting to a subfamily $\mathcal{F}''\subset \mathcal{F}'$ of  $m-1$ vectors among the $m$ Wilczynski vectors $V_{i,j}$, $(i,j)\in\mathcal{F}'$, with the same family $\mathcal{G}$, one has a new reduced Wilczynski matrix, with $m-1$ columns, all of which minors of order $m-1$ are null. So $y_0$ is algebraic relatively to $(\mathcal{F}'',\mathcal{G})$. By induction on  $m$, one has reconstructed  a polynomial $P(x,y)$ vanishing at $y_0$. \\

To obtain the Theorem \ref{theo:wilc}, it suffices now to show that there exists a uniform bound  $N_{d_x,d_y}$ for the depth in  $M_{\mathcal{F},\mathcal{G}}^{red}$ to which we get the reconstruction process, that is, the depth at which we find a first non zero minor. We reach this in the two following lemmas. 

\begin{lemma}\label{lemme:ordreQ}
Let $d_x,\, d_y\in \mathbb{N}^*$. For any series   $y_0=\displaystyle\sum_{n\geq 1} c_nx^{n}\in K[[x]]$ with $c_1\neq 0$, verifying an equation $P(x,y_0)=0$ where $P(x,y)\in K[x,y]$, $\deg_xP\leq d_x,\ \deg_yP\leq d_y$, and for any polynomial $Q(x,y)\in K[x,y],\ \deg_xQ\leq d_x,\ \deg_yQ\leq d_y$, such that $Q(x,y_0)\neq 0$, one has that $\mathrm{ord}_xQ(x,y_0)\leq 2d_xd_y$.
\end{lemma} 
\begin{demo}
Let $y_0$ be a series as in the statement of Lemma \ref{lemme:ordreQ}. We consider the ideal $I_0:=\{R(x,y)\in K[x,y]\ |\ R(x,y_0)=0\}$. By assumption, it is a non trivial prime ideal, so its height is one or two. If it were equal to 2, then it would be a maximal ideal. But $I_0$ is included into the ideal $\{R(x,y)\in K[x,y]\ |\ R(0,0)=0\}$, so:
$$ I_0=\{R(x,y)\in K[x,y]\ |\ R(0,0)=0\}=(x,y)$$
which is absurd because $x\notin I_0$. So, $I_0$ is a height one prime ideal of the factorial ring $K[x,y]$. It is generated by an irreducible polynomial $P_0(x,y)\in K[x,y]$. We set $d_{0,x}:=\deg_x P_0$ and $d_{0,y}:=\deg_y P_0$. Note also that, by factoriality of $K[x,y]$, $P_0$ is also irreducible as an element of $K(x)[y]$.\\
Let $P$ be as in  the statement of Lemma \ref{lemme:ordreQ}. One has that $P=SP_0$ for some $S\in K[x,y]$. Hence $d_{0,x}\leq d_x$ and $d_{0,y}\leq d_y$. Let $Q\in K[x,y]$ such that $Q(x,y_0)\neq 0$ with $\deg_x Q\leq d_x$, $\deg_yQ\leq d_y$. So $P_0$ and $Q$ are coprime in $K(x)[y]$. Their resultant $r(x)$ is non zero. One has the following B\'ezout relation in $K[x][y]$:
$$A(x,y)P_0(x,y)+B(x,y)Q(x,y)=r(x).$$
We evaluate at $y=y_0$:
$$0+B(x,y_0)Q(x,y_0)=r(x).$$
So $\mathrm{ord}_x Q(x,y_0)\leq \deg_x r(x)$. But, the resultant is a  determinant of order at most $d_y+d_{0,y}\leq 2\,d_y$ whose entries are  polynomials in $K[x]$ of degree at most $\max\{d_x,d_{0,x}\}\leq d_x$. So, $ \deg_x r(x)\leq 2\,d_xd_y$. Hence, one has that: $\mathrm{ord}_x Q(x,y_0)\leq  2\,d_xd_y$.
\end{demo}

\begin{lemma}\label{lemme:profondeur}
Let $\mathcal{F}'\subsetneq \mathcal{F}$. If $y_0$ is not algebraic relatively to $(\mathcal{F}',\mathcal{G})$, we denote $l:=|\mathcal{F}'|$ and $p:=\min\left\{k_l\ |\ Q_{\underline{k},\mathcal{F}'}\neq 0,\ \underline{k}=(k_1,\ldots,k_l)\right\}$. Then, for any polynomial $Q(x,y)=\displaystyle\sum_{(i,j)\in\mathcal{F}'\cup\mathcal{G}}b_{i,j}x^iy^j$, we have:
$$\mathrm{ord}_xQ(x,y_0)\leq p\leq 2\,d_xd_y,$$
and the value $p$ is reached for a certain polynomial $Q_0$. 
\end{lemma}
\begin{demo}
By the definition of $p$, for any  $\underline{k}=(k_1,\ldots,k_l)$ with $k_l<p$, we have that $Q_{\underline{k},\mathcal{F}'}=0$. This means that the rank of the column vectors $V_{i,j,p-1}$ that are the restrictions of those of $M_{\mathcal{F}',\mathcal{G}}^{red}$  up to the line $p-1$, is less than $l=|\mathcal{F}'|$. There are coefficients $(a_{i,j})_{(i,j)\in_{\mathcal{F}'\cup\mathcal{G}}}$ not all zero such that $\displaystyle\sum_{(i,j)\in_\mathcal{F}'\cup\mathcal{G}}a_{i,j}V_{i,j,p-1}=(0)$, which is equivalent to the vanishing of the $p-1$ first terms of $Q_0(x,y_0):=\displaystyle\sum_{(i,j)\in\mathcal{F}'\cup\mathcal{G}}a_{i,j}x^i(y_0)^j$. Thus, $\mathrm{ord}_xQ_0(x,y_0)\geq p$, and so $p\leq 2\,d_xd_y$. On the other hand, again by the definition of $p$, the column vectors up to the line $p$ are, in turn, of rank $l=|\mathcal{F}'|$. Any non trivial linear combination is non null, so $\mathrm{ord}_xQ(x,y_0)\leq p$ for all $Q(x,y):=\displaystyle\sum_{(i,j)\in\mathcal{F}'\cup\mathcal{G}}b_{i,j}x^iy^j$.
\end{demo}

We achieve the proof of Theorem \ref{theo:wilc} via the Lemmas  \ref{lemme:ordreQ} and \ref{lemme:profondeur} by considering for a given algebraic series $y_0$ a family  $\mathcal{F}''\subset \mathcal{F}$ minimal among the families such that  $y_0$ is algebraic relatively to $(\mathcal{F}'',\mathcal{G})$. Hence, the  natural number $N_{y_0}$ of (\ref{equ:N}) is always bounded by $N=2\,d_xd_y$.\end{demo}

\noindent\textbf{Construction of the coefficients $a_{i,j}^{(\lambda_0)}$ for a given $y_0$}.\\
Let $y_0$ be algebraic relatively to $(\mathcal{F},\mathcal{G})$ as in \ref{defi:alg-relative}. Let $N=2\,d_xd_y$ as in  \ref{theo:wilc}. We denote by  $M_N$ the matrix consisting in the $N$ first lines of $M_{\mathcal{F},\mathcal{G}}^{red}$. Let $r$ be the rank of $M$, and $m:=r+1$. The minors of $M$ of order $m$ are all zero and there exists a minor of order $m-1=r$ which is non zero. There are two cases. If $r=0$, we choose  $(i,j)\in\mathcal{F}$ and we fix the coefficients $a_{i,j}:=1$ and $a_{l,m}=0$ pour $(l,m)\in\mathcal{F}$, $(l,m)\neq (i,j)$.Then we derive the  coefficients $a_{i,0}$ for $(i,0)\in\mathcal{G}$ from the relations (\ref{equ:terme-cst}). The polynomials $P$ thus obtained are all annihilators of $y_0$. \\
If $r\geq 1$, we consider all the Wilczynski polynomials $Q_{\underline{k},\underline{I}}$ of order $r$ that do not vanish when evaluated at $c_1,\ldots,c_N$. Each of them allows to reconstruct  coefficients $a_{i,j}^{(\lambda)}$, $(i,j)\in\mathcal{F}$, and subsequently coefficients $a_{i,0}$, $(i,0)\in\mathcal{G}$, via (\ref{equ:terme-cst}). The corresponding polynomials $P^{(\lambda)}$ are annihilators of $y_0$ if and only if $\mathrm{ord}_x P^{(\lambda)}\left(x,\, \displaystyle\sum_{k=1}^{N}c_kx^k\right)>N$.

\begin{ex}\label{ex:wilcz2}
We resume the Example \ref{ex:wilc}, and note that, for $\underline{k}=(2,3)$ and for $\underline{I}=((2,1),(2,2))$, we have that:
$$Q_{\underline{k},\underline{I}}=
\left|\begin{array}{cc}c _{1}  & 0 \\c _{2} &  \mathop{\rm c}_{1}^{2} \\
\end{array}\right|=c_1^3\neq 0.$$
So we set $a_{0,2}:=(-1)^2c_1^3=c_1^3$ and, applying  the Cramer's rule :
$$\left\{\begin{array}{lcl}
a_{2,1}&:=&(-1)^1
\left|\begin{array}{cc}2\cdot c _{1}\cdot c _{2}
 & 0 \\ c _{2}^{2}+2~c _{1}~c _{3} & c_{1}^{2} \\ \end{array}\right|=-2c_1^3c_2 \vspace{3pt}\\
 a_{2,2}&:=&(-1)^3 \left|\begin{array}{cc} c _{1} & 2\cdot c _{1}\cdot c _{2}
 \\c _{2} & c _{2}^{2}+2~c _{1}~c _{3} \\ \end{array}\right|=c_1\left(c_2^2-2c_1c_3\right).
\end{array} \right.$$
We deduce from the formulas (\ref{equ:terme-cst}) that:
$$a_{2,0}=-a_{2,1}\cdot 0-a_{0,2}\cdot c_1^2-a_{2,2}\cdot 0=-c_1^5.$$
A vanishing polynomial of a series $y_0=\displaystyle\sum_{n\geq 1} c_nx^{n}\in K((x))$,  $c_1\neq 0$, algebraic relatively to $\mathcal{F}=\left((2,1),(0,2),(2,2)\right)$ and  $\mathcal{G}=(2,0)$ is:
$$\begin{array}{lcl}P(x,y)&=&-c_1^5x^2-2c_1^3c_2\,x^2y+c_1^3\,y^2+c_1\left(c_2^2-2c_1c_3\right)x^2y^2\\
&=&c_1\left[-c_1^4x^2-2c_1^2c_2\,x^2y+c_1^2\,y^2+\left(c_2^2-2c_1c_3\right)x^2y^2\right].\\
\end{array}$$
\end{ex}

\begin{remark}\label{rem:resultant}
\begin{enumerate}
\item Let $y_0$ be a series algebraic with vanishing polynomial of degree $d_x$ in $x$ and $d_y$ in $y$. According to \cite[Chap. 7]{salvy-et-al:cours-mpri}, the method of reconstruction of equation based on Pad\'{e}-Hermite approximants provides a priori only polynomials $P(x,y)=\displaystyle\sum_{i\leq d_x,\ j\leq d_y} a_{i,j}x^iy^j$ such that $P(x,y_0)\equiv 0\, [x^\sigma]$  with $\sigma=(d_x+1)(d_y+1)-1$. Subsequently, one has to check whether $P(x,y_0)=0$ actually. By our Lemma \ref{lemme:ordreQ}, one can always certify that  $P(x,y_0)=0$ just by verifying that $P(x,y_0)\equiv 0\, [x^\tau]$ with $\tau=2d_xd_y$. Hence this reconstruction method as implemented in the GFUN package in Maple software holds for any equation of degree less than $d_x$ in $x$ and $d_y$ in $y$, not for only irreducible ones as in \cite[Theorem 8, p. 110]{salvy-et-al:cours-mpri}.

\item Let us consider the case where $y_0$ is a rational fraction:
$$\begin{array}{lcl}
y_0&=&\displaystyle\frac{-a_0(x)}{a_1(x)}=\displaystyle\frac{-a_{1,0}x-\cdots-a_{d_0,0}x^{d_0}}{1+a_{1,1}x+\cdots+a_{d_1,1}x^{d_1}}\\
&=&\displaystyle\sum_{n\geq 1}c_n x^n\ \textrm{ with }\ c_1\neq 0.
\end{array}$$
Thus, $y_0$ is algebraic relatively to $\mathcal{F}=\left\{(0,1),\ldots,(d_1,1)\right\}$ and $\mathcal{G}=\left\{(1,0),\ldots,\right.$ \\ $\left.(d_0,0)\right\}$. The Wilczynski polynomials of order $|\mathcal{F}|=d_1+1$ are all null. The Wilczynski polynomial $Q_{\underline{k}_0,\underline{I}_0}$ of order $d_1$ with $\underline{k}_0=(1,\ldots,d_1)$ and $\underline{I}_0=((1,1),\ldots,(d_1,1))$ is equal, up to the sign, to the resultant of $a_0(x)$ and $a_1(x)$, by \cite[chap 12 (1.15) p 401]{gelfand-kapranov-zel:discrim-result-det}. 
\item In the present section, the field $K$ can be of any  characteristic. 
\end{enumerate}
\end{remark}


\section{Closed-form expression of an algebraic series.}\label{section:flajo-soria}

Let us assume from now on that $K$ has zero characteristic. Our purpose is to determine the coefficients of an algebraic series in terms of the coefficients of a vanishing polynomial. We consider the following polynomial of degrees bounded by $d_x$ in $x$ and $d_y$ in $y$: 
$$ \begin{array}{lcl}
P(x,y)&=&\displaystyle\sum_{i=0}^{d_x}\displaystyle\sum_{j=0}^{d_y}a_{i,j}x^iy^j ,\ \textrm{  with } P(x,y)\in K[x,y]\\
&=& \displaystyle\sum_{i=0}^{d_x}\pi_{i}(y)x^i\\
&=& \displaystyle\sum_{j=0}^{d_y}a_{j}(x)y^j,
\end{array}$$
and a formal power series which is a simple root:
$$y_0=\displaystyle\sum_{n\geq 1}c_nx^n,  \textrm{  with } y_0\in K[[x]],\   c_1\neq 0.$$
The field $K((x))$ is endowed with the  $x$-adic valuation $\textrm{ord}_x$.

Classically (e.g. \cite{walker_alg-curves}), the resolution of $P=0$ with the Newton-Puiseux method is algorithmic, with two stages:
\begin{enumerate}
\item a first stage of separation of the branches solutions, which illustrates the following fact: $y_0$ may share a principal part with other roots of $P$. This is equivalent to the fact that this principal part is also the principal part of a root of $\partial P/\partial y$.
\item a second stage of unique "automatic" resolution: once the branches are separated, the remaining part of $y_0$ is a root of an  equation called  Henselian in the formal valued context ($y_0$ seen as an algebraic formal power series), and called of  implicit function type in the context of differentiable functions  ($y_0$ seen as the convergent Taylor expansion of  an algebraic function).
\end{enumerate}
We give here a version of the algebraic content of this algorithmic resolution. 

\begin{nota}\label{nota:P_k}
For any $k\in\mathbb{N}$ and for $Q(x,y)=\displaystyle\sum_{j=0}^da_j(x)y^j\in K((x))[y]$, we denote: 
\begin{itemize}
\item $\mathrm{ord}_xQ:=\min\{\mathrm{ord}_x a_j(x),\ j=0,..,d\}$
\item $z_0:=0$ and for $k\geq 1$, $z_k:=\displaystyle\sum_{n=1}^{k}c_nx^n$
\item $y_k:=y_0-z_k=\displaystyle\sum_{n\geq k+1}c_nx^n$
\item $Q_k(x,y):=Q(z_k+x^{k+1}y)=\displaystyle\sum_{i=i_k}^{d_k}\pi_{k,i}(y)x^i$ where $i_k=\mathrm{ord}_x Q_k$ and $d_k:=\deg_xQ_k$
\end{itemize}
\end{nota}

\begin{lemma}\label{lemme:double-simple}
\begin{enumerate}
\item The series $y_0$ is a root of  $P(x,y)$ if and only if the sequence $(i_k)_{k\in\mathbb{N}^*}$ is strictly increasing where $i_k=\mathrm{ord}_x P_k$.
\item The series $y_0$ is a simple root of $P(x,y)$ if and only if the sequence $(i_k)_{k\in\mathbb{N}^*}$ is strictly increasing and there exists a lowest index $k_0$ such that $i_{k_0+1}=i_{k_0}+1$. In that case, one has that  $i_{k+1}=i_k+1$  for any $k\geq k_0$. 
\end{enumerate}
\end{lemma}
\begin{demo}
(1) Note that for any $k$, $i_k\leq\textrm{ord}_xP_k(x,0)=\textrm{ord}_xP(x,z_k)$. Hence, if the sequence $(i_k)_{k\in\mathbb{N}^*}$ is strictly increasing, it tends to  $+\infty$, and so does $\textrm{ord}_xP(x,z_k)$. $y_0$ is indeed a root of $P(x,y)$. Reciprocally, suppose that there exists $1\leq k<l$ such that $i_k\geq i_l$. We apply the Taylor formula to $P_j(x,y)$ for $j>k$:
\begin{center}
\begin{equation}\label{equ:taylor}
\begin{array}{lcl}
P_j(x,y)&=&P_k(x,c_{k+1}+c_{k+2}x+\cdots+x^{j-k}y)\\
&=& \pi_{k,i_k}(c_{k+1})x^{i_k}+\left[\pi_{k,i_k}'(c_{k+1})c_{k+2}+\pi_{k,i_k+1}(c_{k+1})\right]x^{i_k+1}+\cdots.
\end{array}\end{equation}
\end{center}
For $1\leq k\leq j\leq l$, $\, i_l\geq i_j\geq i_k$, so  $i_j=i_k$. Thus,  $\pi_{k,i_k}(c_{k+1})\neq 0$, so for any $j>k$,  $\textrm{ord}_xP_j(x,0)=\textrm{ord}_xP(x,z_j)=i_k$. Hence $\textrm{ord}_xP(x,y_0)=i_k\neq +\infty$.\vspace{2 pt} 

\noindent (2) The series $y_0$ is a double root of $P$ if and only if it is a root of $P$ and $\partial P/\partial y$. We apply the Taylor formula for certain $k\in\mathbb{N}^*$: 
\begin{equation}\label{equ:p_{k+1}}
\begin{array}{lcl}
P_{k+1}(x,y)&=&P_k(x,c_{k+1}+xy)\\
&=& \pi_{k,i_k}(c_{k+1})x^{i_k}+\left[\pi_{k,i_k}'(c_{k+1})y+\pi_{k,i_k+1}(c_{k+1})\right]x^{i_k+1}\\
&& +\left[\displaystyle\frac{\pi_{k,i_{k}}''(c_{k+1})}{2}y^2+\pi_{k,i_{k}+1}'(c_{k+1})\,y+\pi_{k,i_{k}+2}(c_{k+1})\right]x^{i_k+2}+\cdots.
\end{array}
\end{equation}
Note that:
\begin{center}
$\displaystyle\frac{\partial P_k}{\partial y}(x,y)=x^{k+1}\left(\displaystyle\frac{\partial P}{\partial y}\right)_k(x,y)=\displaystyle\sum_{i=i_k}^{d_k}\pi_{k,i}'(y)x^i$
\end{center}
One has that $\pi_{k,i_k}\nequiv 0$ and $\pi_{k,i_k}(c_{k+1})=0$ (see the point (1) above), so $\pi_{k,i_k}'(y)\nequiv 0$. Thus $\mathrm{ord}_x\left(\displaystyle\frac{\partial P}{\partial y}\right)_k=i_k-k-1$. We perform the Taylor expansion of   $\left(\displaystyle\frac{\partial P}{\partial y}\right)_{k+1}=\left(\displaystyle\frac{\partial P}{\partial y}\right)_k(c_{k+1}+x\,y)$:
\begin{center}
$\left(\displaystyle\frac{\partial P}{\partial y}\right)_{k+1}(x,y)=\pi_{k,i_k}'(c_{k+1})x^{i_k-k-1}+\left[\pi_{k,i_k}''(c_{k+1})y+\pi_{k,i_k+1}' (c_{k+1})\right]x^{i_k-k}+\cdots$.
\end{center}
By the point (1), if $y_0$ is a double root $P$, we must have $\pi_{k,i_k}'(c_{k+1})=0$. Moreover, if $\pi_{k,i_k+1}(c_{k+1})\neq 0$, we would have $i_{k+1}=i_k+1$ and even $i_{k+j}=i_k+1$ for every  $j$ according to (\ref{equ:taylor}): $y_0$ could not be a root of $P$. So, $\pi_{k,i_k+1}(c_{k+1})= 0$, and, accordingly, $i_{k+1}\geq i_k+2$.\\
If $y_0$ is a simple root of $P$, from the point (1) there exists a lowest  natural number $k_0$ such that the sequence $(i_k-k-1)_{k\in\mathbb{N}^*}$ is no longer strictly increasing, that is, such that  $\pi_{k_0,i_{k_0}}'(c_{k_0+1})\neq 0$. For any $k\geq k_0$, we consider the  Taylor expansion of $\left(\displaystyle\frac{\partial P}{\partial y}\right)_{k+1}=\left(\displaystyle\frac{\partial P}{\partial y}\right)_{k_0}(c_{k_0+1}+\cdots+x^{k-k_0+1}y)$:
\begin{equation}\label{equ:partialP}
\left(\displaystyle\frac{\partial P}{\partial y}\right)_{k+1}(x,y)=\pi_{k_0,i_{k_0}}'(c_{k_0+1})x^{i_{k_0}-k_0-1}+\left[\pi_{k_0,i_{k_0}}''(c_{k_0+1})c_{k_0+2}+\pi_{k_0,i_{k_0}+1}' (c_{k_0+1})\right]x^{i_{k_0}-k_0}+\cdots
\end{equation}
and we get that:
\begin{equation}\label{equ:taylor-deriv}
\mathrm{ord}_x\left(\displaystyle\frac{\partial P}{\partial y}\right)_{k+1}(x,0)=\mathrm{ord}_x\left(\displaystyle\frac{\partial P}{\partial y}\right)_{k+1}=i_{k_0}-k_0-1.
\end{equation} As $\pi_{k+1,i_{k+1}}'(y)\nequiv 0$, we obtain that $i_{k+1}=\mathrm{ord}_xP_{k+1}=\mathrm{ord}_x\left(\displaystyle\frac{\partial P_{k+1}}{\partial y}\right)=k+2+\mathrm{ord}_x\left(\displaystyle\frac{\partial P}{\partial y}\right)_{k+1}=i_{k_0}+k-k_0+1$. Hence, from the rank $k_0$, the sequence $(i_k)$ increases  one by one. 
\end{demo}

Resuming the notations of the Theorem \ref{theo:wilc} and of the Lemma \ref{lemme:double-simple}, the  natural number $k_0$ represents the length of the principal part in the stage of separation of the branches. In the following lemma, we bound it using the Lemma \ref{lemme:ordreQ} or the discriminant $\Delta_P$ of $P$. 

\begin{lemma}\label{lemme:partie-princ}
With the notations of the  Theorem \ref{theo:wilc}, the  natural number $k_0$ verifies that:
$$k_0\leq 2d_x d_y+1.$$
In particular, if $P$ has only simple roots: 
$$k_0\leq d_x(2\,d_y-1)+1.$$
 \end{lemma}
 \begin{demo}
By the Lemma \ref{lemme:ordreQ}, since $P(x,y_0)=0$ and $ \displaystyle\frac{\partial P}{\partial y}(x,y_0)\neq 0$, one has that:
$$\mathrm{ord}_x \displaystyle\frac{\partial P}{\partial y}(x,y_0)\leq 2d_xd_y.$$
But, by definition, $k_0$ is the lowest  natural number such that: $$\mathrm{ord}_x \displaystyle\frac{\partial P}{\partial y}(x,z_{k_0+1})=\mathrm{ord}_x\displaystyle\frac{\partial P}{\partial y}(x,z_{k_0})=i_{k_0}-k_0-1$$ (see the point (2) of the preceding lemma). Hence, we get from (\ref{equ:taylor-deriv}) that $\mathrm{ord}_x \displaystyle\frac{\partial P}{\partial y}(x,z_{k})=\mathrm{ord}_x\displaystyle\frac{\partial P}{\partial y}(x,z_{k_0})=i_{k_0}-k_0-1$  for any $k> k_0$. So, $\mathrm{ord}_x \displaystyle\frac{\partial P}{\partial y}(x,y_0)=\mathrm{ord}_x\displaystyle\frac{\partial P}{\partial y}(x,z_{k_0})$. To conclude, we note that $\mathrm{ord}_x\displaystyle\frac{\partial P}{\partial y}(x,z_{k_0})\geq k_0-1$ by definition of $k_0$.\\
In the case where  $P$  has only  simple roots, as in the proof of  the Lemma \ref{lemme:ordreQ},  $\mathrm{ord}_x \displaystyle\frac{\partial P}{\partial y}(x,y_0)$ is bounded by the degree of the resultant of $P$ and $\displaystyle\frac{\partial P}{\partial y}$, say the discriminant $\Delta_P$  of $P$, which is bounded by $d_x(2\,d_y-1)+1$. 
 \end{demo}

\begin{theo}\label{theo:FS}
Consider the following polynomial in $K[x,y]$ of given degrees $d_x$ in $x$ and $d_y$ in $y$: 
$$ P(x,y)=\displaystyle\sum_{i=0}^{d_x}\displaystyle\sum_{j=0}^{d_y}a_{i,j}x^iy^j = \displaystyle\sum_{i=0}^{d_x}\pi_{i}(y)x^i,$$
and a formal power series which is a simple root:
$$y_0=\displaystyle\sum_{n\geq 1}c_nx^n\ \in K[[x]],\   c_1\neq 0.$$
Resuming the  notations of \ref{nota:P_k} and \ref{lemme:double-simple}, we set $\omega_0:=\pi_{k_0,i_{k_0}}'(c_{k_0+1})\neq 0$.
Hence, for any $k>k_0$ :
\begin{itemize}
\item either the polynomial $z_{k+1}=\displaystyle\sum_{n=1}^{k+1}c_nx^n$ is a solution of $P(x,y)=0$;
\item or the polynomial $R_k(x,y):=\displaystyle\frac{P_k(x,y+c_{k+1})}{-\omega_0x^{i_k}}=-y+Q_k(x,y)$  defines a reduced Henselian equation:
$$ y=Q_k(x,y)$$
with $Q_k(0,y)\equiv 0$ and satisfied by:
$$t_{k+1}:=\frac{y_0-z_{k+1}}{x^{k+1}}=c_{k+2}x+c_{k+3}x^2+\cdots.$$
\end{itemize} 
\end{theo}
\begin{demo}
We show by induction on  $k>k_0$ that $R_k(x,y)=-y+xT_k(x,y)$ with $T_k(x,y)\in K[x,y]$. For  $k=k_0+1$, by (\ref{equ:p_{k+1}}), since $i_{k_0+1}=i_{k_0}+1$, we have that:
$$P_{k_0+1}(x,y)=\left[\omega_0\,y+\pi_{k_0,i_{k_0}+1}(c_{k_0+1})\right]x^{i_{k_0}+1}+\cdots.$$
Since $i_{k_0+2}=i_{k_0}+2$, $\pi_{k_0+1,i_{k_0}+1}(y)=\omega_0\,y+\pi_{k_0,i_{k_0}+1}(c_{k_0+1})$ vanishes at $c_{k_0+2}$, which implies that $c_{k_0+2}=\displaystyle\frac{-\pi_{k_0,i_{k_0}+1}(c_{k_0+1})}{\omega_0}$. Computing $R_{k_0+1}(x,y)$, it follows that:
\begin{center}
$R_{k_0+1}(x,y)=-y+Q_{k_0+1}(x,y) \textrm{ with }$\end{center}
$$ Q_{k_0+1}(x,y)=x\left[\displaystyle\frac{\pi_{k,i_{k_0}}''(c_{k_0+1})}{2}(y+c_{k_0+2})^2+\pi_{k_0,i_{k_0}+1}'(c_{k_0+1})\,(y+c_{k_0+2})+\pi_{k_0,i_{k_0}+2}(c_{k_0+1})\right]+x^2[\cdots].$$
So $Q_{k_0+1}(0,y)\equiv 0$. 

Suppose that the property holds true at a rank $k>k_0+1$. It follows that:
$$\begin{array}{lcl}
 P_k(x,y)&=&\omega_0(y-c_{k+1})x^{i_k}+x^{i_k+1}\tilde{T}_k(x,y)\\
 &=&\pi_{k,i_k}(y)x^{i_k}+\pi_{k,i_k+1}(y)x^{i_k+1}+ \cdots.
\end{array}$$
Since $P_{k+1}(x,y)=P_k(x,c_{k+1}+xy)$, we have that:
$$ P_{k+1}(x,y)=\left[\omega_0\,y+\pi_{k,i_k+1}(c_{k+1})\right]x^{i_k+1}+\pi_{k+1,i_k+2}(y)x^{i_k+2}+\cdots.$$
But $i_k+2=i_{k+2}>i_{k+1}=i_k+1$. So we must have $\pi_{k+1,i_k+1}(c_{k+2})=0$. So, $c_{k+2}=\displaystyle\frac{-\pi_{k,i_k+1}(c_{k+1})}{\omega_0}$. It follows that:
$$P_{k+1}(x,y)=\omega_0(y-c_{k+2})x^{i_k+1}+\pi_{k+1,i_k+2}(y)x^{i_k+2}+\cdots,$$
Hence:
$$\begin{array}{lcl}
R_{k+1}(x,y)&=&-y-x\displaystyle\frac{\pi_{k+1,i_k+2}(y+c_{k+2})}{\omega_0}+x^2[\cdots]+\cdots\\
&=&-y+xT_{k+1}(x,y),\ \ \ \ T_{k+1}\in K[x,y],
\end{array} $$
as desired.

In particular, $Q_{k}(0,0)=\displaystyle\frac{\partial Q_{k}}{\partial y}(0,0)=0$. So the equation $y=Q_k(x,y)$ is reduced Henselian if and only if $Q_k(x,0)\nequiv 0$, which is equivalent to $z_{k+1}$ not being a root of $P$.
\end{demo}

\begin{remark}
By (\ref{equ:partialP}), we note that $$\left(\displaystyle\frac{\partial P}{\partial y}\right)(x,y_0)=\omega_0x^{i_{k_0}-k_0-1}+\cdots.$$
Thus, $\omega_0$ is the  initial coefficient of $\left(\displaystyle\frac{\partial P}{\partial y}\right)(x,y_0)$.
 \end{remark}

For the courageous reader, in the case where $y_0$ is a series which is not a polynomial, we deduce from \ref{theo:FS} and from the   Flajolet-Soria formula \ref{theo:formule-FS} a closed-form expression for the coefficients of $y_0$ in terms of the coefficients $a_{i,j}$ of $P$ and of the coefficients of an initial part  $z_k$ of $y_0$ sufficiently large.

\begin{coro}\label{coro:FS}
For any $k\geq k_0+1$, for any $p\geq 1$, one has that:
$$c_{k+1+p}=\displaystyle\sum_{q=1}^p\displaystyle\frac{1}{q}\left(\displaystyle\frac{-1}{\omega_0}\right)^q\displaystyle\sum_{|S|=q,\ \|S\|_2\geq q-1}A^S\left(\displaystyle\sum_{|T_S|=\|S\|_2-q+1 \atop \|T_S\|=p+qi_k-(q-1)(k+1)-\|S\|_1}e_{T_S}C^{T_S}\right),$$
where $S=(s_{i,j})$, $A^S=\displaystyle\prod_{i=0,\ldots,d_x\, ,\ j=0,\ldots,d_y}a_{i,j}^{s_{i,j}}$, $T_S=(t_{S,i})$,  $C^{T_S}=\displaystyle\prod_{i=1}^{k+1}c_i^{t_{S,i}}$, 
 and $e_{T_S}\in\mathbb{N}$ is of the form:
$$e_{T_S}=\displaystyle\sum_{n^{l,m}_{i,j,L}}\displaystyle\frac{q!}{\displaystyle\prod_{l=1,\ldots,(k+1)d_y+d_x-i_k \atop m=0,\ldots,m_l}\displaystyle\prod_{i=0,\ldots,d_x \atop j=m,\ldots,d_y}\displaystyle\prod_{|L|=j-m \atop \|L\|=l+i_k-m(k+1)-i}n^{l,m}_{i,j,L}!}\displaystyle\prod_{l=1,\ldots,(k+1)d_y+d_x-i_k \atop m=0,\ldots,m_l}\displaystyle\prod_{i=0,\ldots,d_x \atop j=m,\ldots,d_y}\displaystyle\prod_{|L|=j-m \atop \|L\|=l+i_k-m(k+1)-i}\left(\displaystyle\frac{j!}{m!\,L!}\right)^{n^{l,m}_{i,j,L}},$$
where $m_l:=\min\left\{\left\lfloor\displaystyle\frac{l+i_k}{k+1}\right\rfloor,d_y\right\}$, $L=L_{i,j}^{l,m}=\left(l_{i,j,1}^{l,m},\ldots,l_{i,j,k+1}^{l,m}\right)$, 
 and where the sum is taken over the set of  $\left(n^{l,m}_{i,j,L}\right)_{ l=1,\ldots,(k+1)d_y+d_x-i_k,\ m=0,\ldots,m_l \atop i=0,\ldots,d_x,\ j=m,\ldots,d_y,\ |L|=j-m,\  \|L\|=l+i_k-m(k+1)-i}$ such that:
\begin{center}
 $\displaystyle\sum_{l,m}\displaystyle\sum_{i,j}\displaystyle\sum_{L}n^{l,m}_{i,j,L}=q\ \ \ $ and $\ \ \ \displaystyle\sum_{l,m}\displaystyle\sum_{i,j}\displaystyle\sum_{L}n^{l,m}_{i,j,L}L=T_S$.
\end{center}

\end{coro}
\begin{demo}
We get started by computing the coefficients of $\omega_0x^{i_k}R_k$,  in order to get those of $Q_k$:
$$\begin{array}{lcl}
\omega_0x^{i_k}R_k&=&P_k(x,\,y+c_{k+1})\\
&=&P(x,z_{k+1}+x^{k+1}y)\\
&=& \displaystyle\sum_{i=0,\ldots,d_x\, ,\ j=0,\ldots,d_y}a_{i,j}x^i\left(z_{k+1}+x^{k+1}y\right)^{j}\\
&=& \displaystyle\sum_{i=0,\ldots,d_x\, ,\ j=0,\ldots,d_y}a_{i,j}x^i\displaystyle\sum_{m=0}^{j}\displaystyle\frac{j!}{m!\,(j-m)!}z_{k+1}^{j-m}x^{m(k+1)}y^m.
\end{array}$$
For $L=(l_1,\cdots,l_{k+1})$, we denote 
$C^L:=c_1^{l_1}\cdots c_{k+1}^{l_{k+1}}$. One  has that:
$$z_{k+1}^{j-m}=\displaystyle\sum_{|L|=j-m}\displaystyle\frac{(j-m)!}{L!}C^Lx^{\|L\|}.$$
So:
$$ \omega_0x^{i_k}R_k=\displaystyle\sum_{m=0}^{d_y} \displaystyle\sum_{i=0,\ldots,d_x \atop j=m,\ldots,d_y}a_{i,j}\displaystyle\sum_{|L|=j-m}\displaystyle\frac{j!}{m!\,L!}C^Lx^{\|L\|+m(k+1)+i}y^m.$$
We set $\hat{l}=\|L\|+m(k+1)+i$, which ranges between $m(k+1)$ and $(k+1)(d_y-m)+m(k+1)+d_x=(k+1)d_y+d_x$. Thus: 
$$\omega_0x^{i_k}R_k=\displaystyle\sum_{m=0,\ldots,d_y \atop \hat{l}=m(k+1),\ldots,(k+1)d_y+d_x} \displaystyle\sum_{i=0,\ldots,d_x \atop j=m,\ldots,d_y}a_{i,j}\displaystyle\sum_{|L|=j-m \atop \|L\|=\hat{l}-m(k+1)-i}\displaystyle\frac{j!}{m!\,L!}C^Lx^{\hat{l}}y^m.$$
Since $R_k=-y+Q_k(x,y)$ with $Q_k(0,y)\equiv 0$, the coefficients of $Q_k$ are obtained for $\hat{l}=i_k+1,\ldots,(k+1)d_y+d_x$. We set $l:=\hat{l}-i_k$, $m_l:=\min\left\{\left\lfloor\displaystyle\frac{l+i_k}{k+1}\right\rfloor,d_y\right\}$ and we have that: $$Q_k(x,y)=\displaystyle\sum_{l=1,\ldots,(k+1)d_y+d_x-i_k \atop m=0,\ldots,m_l}b_{l,m}x^ly^m,$$
with:
$$b_{l,m}=\displaystyle\frac{-1}{\omega_0}\displaystyle\sum_{i=0,\ldots,d_x \atop j=m,\ldots,d_y}a_{i,j}\displaystyle\sum_{|L|=j-m \atop \|L\|=l+i_k-m(k+1)-i}\displaystyle\frac{j!}{m!\,L!}C^L.$$
We are in position to apply  the version \ref{rem:FS}  of the Flajolet-Soria formula \ref{theo:formule-FS} in order to compute the coefficients of $t_k=c_{k+2}x+c_{k+3}x^2+\cdots$.  Thus, denoting $Q:=(q_{l,m})$ for $l=1,\ldots,(k+1)d_y+d_x-i_k$ and $m=0,\ldots,m_l$, we get that: 
$$c_{k+1+p}=\displaystyle\sum_{q=1}^p\displaystyle\frac{1}{q}\displaystyle\sum_{|Q|=q,\, \|Q\|_1=p,\,\|Q\|_2=q-1}\displaystyle\frac{q!}{Q!}B^Q.$$
Let us compute:
$$\begin{array}{lcl}
b_{l,m}^{q_{l,m}}&=&\left(\displaystyle\frac{-1}{\omega_0}\right)^{q_{l,m}}\left(\displaystyle\sum_{i=0,\ldots,d_x \atop j=m,\ldots,d_y}a_{i,j}\displaystyle\sum_{|L|=j-m \atop \|L\|=l+i_k-m(k+1)-i}\displaystyle\frac{j!}{m!\,L!}C^L\right)^{q_{l,m}}\\
&=&\left(\displaystyle\frac{-1}{\omega_0}\right)^{q_{l,m}}\displaystyle\sum_{|M_{l,m}|=q_{l,m},\, \|M_{l,m}\|_2\geq m\,q_{l,m}}\displaystyle\frac{q_{l,m}!}{M_{l,m}!} A^{M_{l,m}}\displaystyle\prod_{i=0,\ldots,d_x \atop j=m,\ldots,d_y}\left(\displaystyle\sum_{|L|=j-m \atop \|L\|=l+i_k-m(k+1)-i}\displaystyle\frac{j!}{m!\,L!}C^L\right)^{m^{l,m}_{i,j}}\\
&& \textrm{ where } M_{l,m}=(m^{l,m}_{i,j})\textrm{ for } i=0,\ldots,d_x,\  j=m,\ldots,d_y.
\end{array}
$$ 
For each $m^{l,m}_{i,j}$, we enumerate the terms $\displaystyle\frac{j!}{m!\,L!}C^L$ with $u=1,\ldots,\alpha_{i,j}$. Subsequently:
$$\begin{array}{lcl}
\left(\displaystyle\sum_{|L|=j-m \atop \|L\|=l+i_k-m(k+1)-i}\displaystyle\frac{j!}{m!\,L!}C^L\right)^{m^{l,m}_{i,j}}&=& \left(\displaystyle\sum_{u=1}^{\alpha_{i,j}}\displaystyle\frac{j!}{m!\,L_u!}C^{L_u}\right)^{m^{l,m}_{i,j}}\\
&=& \displaystyle\sum_{|N_{i,j}|=m^{l,m}_{i,j}}\displaystyle\frac{m^{l,m}_{i,j}!}{N_{i,j}!}\displaystyle\prod_{u=1}^{\alpha_{i,j}}\left(\displaystyle\frac{j!}{m!\,L_u!}\right)^{n^{l,m}_{i,j,u}}C^{\sum_{u=1}^{\alpha_{i,j}}n^{l,m}_{i,j,u}L_u},
\end{array}
 $$
where $N^{l,m}_{i,j}=\left(n^{l,m}_{i,j,u}\right)_{u=1,\ldots,\alpha_{i,j}}$. 
 Denoting $ U_{l,m}:=\displaystyle\sum_{i=0,\ldots,d_x \atop j=m,\ldots,d_y}\displaystyle\sum_{u=1}^{\alpha_{i,j}}n^{l,m}_{i,j,u}L_u$, one has that:
\begin{center}$\begin{array}{lcl}
|U_{l,m}|&=&\displaystyle\sum_{i=0,\ldots,d_x \atop j=m,\ldots,d_y}\displaystyle\sum_{u=1}^{\alpha_{i,j}}n^{l,m}_{i,j,u}|L_u|\\
&=&\displaystyle\sum_{i=0,\ldots,d_x \atop j=m,\ldots,d_y}\left(\displaystyle\sum_{u=1}^{\alpha_{i,j}}n^{l,m}_{i,j,u}\right)(j-m)\\
&=&\displaystyle\sum_{i=0,\ldots,d_x \atop j=m,\ldots,d_y}m^{l,m}_{i,j}(j-m)\\
&=& \|M_{l,m}\|_2-m\,q_{l,m}.
\end{array}$\end{center} 
Likewise, one has that:
\begin{center}$\begin{array}{lcl}
\|U_{l,m}\|&=&\displaystyle\sum_{i=0,\ldots,d_x \atop j=m,\ldots,d_y}\displaystyle\sum_{u=1}^{\alpha_{i,j}}n^{l,m}_{i,j,u}\|L_u\|\\
&=&\displaystyle\sum_{i=0,\ldots,d_x \atop j=m,\ldots,d_y}\left(\displaystyle\sum_{u=1}^{\alpha_{i,j}}n^{l,m}_{i,j,u}\right)(l+i_k-m(k+1)-i)\\
&=&\displaystyle\sum_{i=0,\ldots,d_x \atop j=m,\ldots,d_y}m^{l,m}_{i,j}(l+i_k-m(k+1)-i)\\
&=& q_{l,m}[l+i_k-m(k+1)]-\|M_{l,m}\|_1.
\end{array}$\end{center} 
So we get that:
\begin{center}
$\begin{array}{lcl}
b_{l,m}^{q_{l,m}}&=& \left(\displaystyle\frac{-1}{\omega_0}\right)^{q_{l,m}}\displaystyle\sum_{|M_{l,m}|=q_{l,m},\, \|M_{l,m}\|_2\geq m\,q_{l,m}}A^{M_{l,m}}\displaystyle\sum_{|U_{l,m}|=\|M_{l,m}\|_2-m\,q_{l,m} \atop \|U_{l,m}\|=q_{l,m}[l+i_k-m(k+1)]-\|M_{l,m}\|_1}d_{U_{l,m}}C^{U_{l,m}}\\
&& \textrm{ with }d_{U_{l,m}}:=\displaystyle\sum_{N^{l,m}_{i,j}}\displaystyle\frac{q_{l,m}!}{\displaystyle\prod_{i=0,\ldots,d_x \atop j=m,\ldots,d_y}N^{l,m}_{i,j}!}\displaystyle\prod_{i=0,\ldots,d_x \atop j=m,\ldots,d_y}\displaystyle\prod_{u=1}^{\alpha_{i,j}}\left(\displaystyle\frac{j!}{m!\,L_u!}\right)^{n^{l,m}_{i,j,u}},
 \end{array}$
 \end{center} 
 \textrm{ where the sum is taken over }$\left\{\left(N^{l,m}_{i,j}\right)_{i=0,\ldots,d_x \atop j=m,\ldots,d_y}\textrm{ such that }|N^{l,m}_{i,j}|=m^{l,m}_{i,j}\textrm{ and }\displaystyle\sum_{i=0,\ldots,d_x \atop j=m,\ldots,d_y}\sum_{u=1}^{\alpha_{i,j}}n^{l,m}_{i,j,u}L_u=U_{l,m}\right\}$.

We deduce that:
$$\begin{array}{lcl}
B^{Q}&=&\displaystyle\prod_{l=1,\ldots,(k+1)d_y+d_x-i_k,\, m=0,\ldots,m_l}b_{l,m}^{q_{l,m}}\\
&=& \left(\displaystyle\frac{-1}{\omega_0}\right)^{q}\displaystyle\prod_{l,m}\left[\displaystyle\sum_{|M_{l,m}|=q_{l,m},\, \|M_{l,m}\|_2\geq m\,q_{l,m}}A^{M_{l,m}}\displaystyle\sum_{|U_{l,m}|=\|M_{l,m}\|_2-m\,q_{l,m} \atop \|U_{l,m}\|=q_{l,m}[l+i_k-m(k+1)]-\|M_{l,m}\|_1}d_{U_{l,m}}C^{U_{l,m}}\right].
\end{array}
$$ 
We set  $S:=\displaystyle\sum_{l,m}M_{l,m}$. So $|S|=\displaystyle\sum_{l,m}q_{l,m}=q$ and $\|S\|_2\geq \displaystyle\sum_{l,m}mq_{l,m}=\|Q\|_2=q-1$. Moreover, for any fixed $S$, we set $T_S:=\displaystyle\sum_{l,m}U_{l,m}$. So $
|T_S|=\displaystyle\sum_{l,m} \|M_{l,m}\|_2-m\,q_{l,m} =\|S\|_2-\|Q\|_2=\|S\|_2-q+1$, and:
\begin{center}  $\begin{array}{lcl}
\|T_S\|&=&\displaystyle\sum_{l,m}q_{l,m}[l+i_k-m(k+1)]-\|M_{l,m}\|_1\\
&=&\|Q\|_1+|Q|i_k-\|Q\|_2(k+1)-\|S\|_1\\
&=&p+qi_k-(q-1)(k+1)-\|S\|_1.
\end{array}$\end{center} Thus, as desired:
$$\begin{array}{lcl}
\displaystyle\sum_{|Q|=q,\, \|Q\|_1=p,\,\|Q\|_2=q-1}\displaystyle\frac{q!}{Q!}B^{Q}&=&  \left(\displaystyle\frac{-1}{\omega_0}\right)^{q}\displaystyle\sum_{|S|=q,\, \|S\|_2\geq q-1}A^{S}\displaystyle\sum_{|T_S|=\|S\|_2-q+1 \atop \|T_S\|=p+qi_k-(q-1)(k+1)-\|S\|_1}e_{T_S}C^{T_S},\end{array}
$$ 
\textrm{ where } $e_{T_S}:=\displaystyle\sum_{N^{l,m}_{i,j}}\displaystyle\frac{q!}{\displaystyle\prod_{l,m}\displaystyle\prod_{i,j}N^{l,m}_{i,j}!}\displaystyle\prod_{l,m} \displaystyle\prod_{i,j}\displaystyle\prod_{u}\left(\displaystyle\frac{j!}{m!\,L_u!}\right)^{n^{l,m}_{i,j,u}}$  and  \textrm{ where the sum is taken over }\\ $\left\{\left(N^{l,m}_{i,j}\right)_{ l=1,\ldots,(k+1)d_y+d_x-i_k,\, m=0,\ldots,m_l \atop i=0,\ldots,d_x,\ j=m,\ldots,d_y}\textrm{ such that }\displaystyle\sum_{l,m}\displaystyle\sum_{i,j}|N^{l,m}_{i,j}|=q\textrm{ and } \displaystyle\sum_{l,m}\displaystyle\sum_{i,j} \sum_{u=1}^{\alpha_{i,j}}n^{l,m}_{i,j,u}L_u=T_S\right\}$.

\end{demo}

\begin{remark}\label{rem:omega_0}
We have seen in the  Theorem \ref{theo:FS} and its proof that $\omega_0=\pi'_{k_0,i_{k_0}}(c_{k_0+1})$ is the coefficient of the monomial $x^{i_{k_0}+1}y$ in the expansion of $P_{k_0+1}(x,y)=P(x,c_1x+\cdots+c_{k_0+1}x^{k_0+1}+x^{k_0+2}y)$, and that  $c_{k_0+2}=\displaystyle\frac{-\pi_{k_0,i_{k_0}+1}(c_{k_0+1})}{\omega_0}$ where $\pi_{k_0,i_{k_0}+1}(c_{k_0+1})$ is the coefficient of $x^{i_{k_0}+1}$ in the expansion of $P_{k_0+1}(x,y)$. Expanding $P_{k_0+1}(x,y)$, having done the whole computations, we deduce that:
$$\left\{\begin{array}{lcl}
\omega_0&=&\displaystyle\sum_{i=0,..,d_x,\ j=0,..,d_y}\ \ \displaystyle\sum_{|L|=j,\  \|L\|=i_{k_0}+1-i}\displaystyle\frac{j!}{L!}a_{i,j}C^L\ ;\\
c_{k_0+2}&=& \displaystyle\frac{-1}{\omega_0}\displaystyle\sum_{i=0,..,d_x,\ l=0,..,d_y-1}\ \ \displaystyle\sum_{|L|=l,\  \|L\|=i_{k_0}+k_0-i-1}\ \ \displaystyle\frac{(l+1)!}{L!}a_{i,l+1}C^L.
\end{array}\right. $$
\end{remark}

\begin{ex}\label{ex:FS}
In order to illustrate the Corollary \ref{coro:FS} and its proof, we resume the polynomial of the Example \ref{ex:wilc}:
$$\begin{array}{lcl}
P(x,y)&=&a_{0,2}y^2+\left( a_{2,0} +a_{2,1}y+a_{2,2}y^2\right)x^2\\
P_0(x,y)&=&\left(a_{2,0}+a_{0,2}y^2\right)x^2+a_{2,1}yx^3+a_{2,2}y^2x^4\\
P_1(x,y)&=&\left( 2a_{0,2}c_1y+a_{2,1}c_1\right)x^3+\left(a_{0,2}y^2+a_{2,1}y+a_{2,2}c_1^2\right)x^4 +2a_{2,2}c_1yx^5+a_{2,2}y^2x^6\\
&&\textrm{with }a_{2,0}+a_{0,2}c_1^2=0\,\Leftrightarrow\, c_1=\pm\sqrt{\displaystyle\frac{-a_{2,0}}{a_{0,2}}}.
\end{array}$$
 Thus, $i_0=2$, $i_1=3=i_0+1$, so $k_0=0$, $\omega_0=2\,a_{0,2}\,c_1$. The coefficient $c_2$ must verify $2a_{0,2}c_1c_2+a_{2,1}c_1=0\, \Leftrightarrow\, c_2=\displaystyle\frac{-a_{2,1}}{2a_{0,2}}$. We obtain that:
\begin{center}
$\omega_0\,R_1= \omega_0y+  \left(a_{{2,2}}{c_{{1}}}^{2}+ a_{{2,1}}c_{{2}}+a_{{0,2}}{c_{{2
}}}^{2}+\left(a_{{2,1}}+2\,a_{{0,2}}c_{{2}}\right)y+a_{{0,2}}{y}^{2}
 \right) {x}+ \left( 2\,a_{{2,2}}c_{{1}}c_{{2}}+2\,a_{{2,2}}c_{{1}
}y \right) {x}^{2}+ \left( a_{{2,2}}{c_{{2}}}^{2}+2\,a_{{2,2}}c_{{2}}y+a_{{2,2}}{y}^{2}
 \right) {x}^{3}.$
 \end{center}
So the coefficients of the corresponding reduced Henselian equation $y=Q_1(x,y)$ are:
 \begin{center}
 $\begin{array}{ll}
 b_{1,0}=-\left(a_{{2,2}}{c_{{1}}}^{2}+ a_{{2,1}}c_{{2}}+a_{{0,2}}{c_{{2
}}}^{2}\right)/\omega_0, & b_{1,1}=-\left(a_{{2,1}}+2\,a_{{0,2}}c_{{2}}\right)/\omega_0=0,
 \end{array}$
 
$ \begin{array}{lll}
 b_{1,2}=-a_{{0,2}}/\omega_0,& b_{2,0}= -2\,a_{{2,2}}c_{{1}}c_{{2}}/\omega_0,& b_{2,1}=-2\,a_{{2,2}}c_{{1}}/\omega_0,\\
 b_{3,0}=-a_{{2,2}}{c_{{2}}}^{2}/\omega_0, &b_{3,1}=-2\,a_{{2,2}}c_{{2}}/\omega_0, &b_{3,2}=-a_{{2,2}}/\omega_0,
 \end{array}$\end{center}
But, by the version \ref{rem:FS} of the Flajolet-Soria formula \ref{theo:formule-FS}, one has that:\\
$\begin{array}{lcl}
c_3&=&b_{1,0}=\displaystyle\frac{- a_{{2,2}}{c_{{1}}}^{2}- a_{{2,1}}c_{{2}}-a_{{0,2}}{c_{{2
}}}^{2} }{2\,a_{0,2}\,c_1};\\
c_4&=&b_{2,0}+b_{1,0}b_{1,1}=b_{2,0}=\displaystyle\frac{-2\,a_{{2,2}}c_{{1}}c_{{2}}}{2\,a_{0,2}\,c_1};\\
c_5&=&b_{3,0}+b_{1,0}b_{2,1}+b_{1,0}^2b_{1,2}+b_{1,0}b_{1,1}^2+b_{2,0}b_{1,1}=b_{3,0}+b_{1,0}b_{2,1} +b_{1,0}^2b_{1,2}\\
&=&\displaystyle\frac{-a_{2,2}\,{c_2}^2}{2\,a_{0,2}\,c_1}+ \displaystyle\frac{2\,a_{{2,1}}a_{{2,2}}c_{{1}}c_{{2}}+2\,a_{{0,2}}a_{{2,2}}c_{{1}}{c_{{2
}}}^{2}+2\,{a_{{2,2}}}^{2}{c_{{1}}}^{3}
}{\left(2\,a_{0,2}\,c_1\right)^2}\,-
\end{array}\\
\begin{array}{l} \displaystyle\frac{a_{{0,2}}{a_{{2,1}}}^{2}{c_{{2}}}^{2}+2\,{a_{{0,2}}}^{2}a_{{2,1}}{c_{{2}}}^{3}+2\,a_{{0,2}}a_{{2,1}}a_{{2,2}}{c_{{1}}}^{2}c_{{2}}+{a_{{0,2}}
}^{3}{c_{{2}}}^{4}+2\,{a_{{0,2}}}^{2}a_{{2,2}}{c_{{1}}}^{2}{c_{{2}}}^{
2}+a_{{0,2}}{a_{{2,2}}}^{2}{c_{{1}}}^{4} }{\left(2\,a_{0,2}\,c_1\right)^3};\\
\ \ \ \vdots\end{array}$
\end{ex}

\begin{remark}\label{rem:furstenberg}
Classically, a series $y_0=\displaystyle\sum_{n\geq 0} c_nx^n\in K[[x]]$ is algebraic if and only if its coefficients $c_n$ are the diagonal coefficients of the power series expansion of a bivariate rational fraction \cite{furstenberg:alg-funct,denef-lipshitz:alg-series-diag}. In particular, in the reduced Henselian case $y=Q(x,y)$ (see \ref{defi:equ-hensel-red}), the rational fraction can be written:
$$y_0=\mathrm{Diag}\left( \displaystyle\frac{y^2-y^2\displaystyle\frac{\partial Q}{\partial y}(xy,y)}{y-Q(xy,y)}\right).$$
With the computations  in the proof of the Corollary \ref{coro:FS}, we can deduce in the general case $P(x,y)=0$ a formula for the  rational fraction having the $c_n$ as diagonal coefficients of its expansion.

\end{remark}

As a consequence of the Theorem \ref{theo:wilc} and Corollary \ref{coro:FS}, we get the following result:

\begin{coro}
Let $d_x,d_y\in \mathbb{N}^*$. We set $\mu:=2d_xd_y+2$ and $M:=\displaystyle\frac{1}{2}d_y(d_y+1)(d_x+1)+d_x+d_y$. There exists a finite set $\Lambda$ and for any $\lambda\in\Lambda$, there exist a polynomial $\Omega^{(\lambda)}(x_1,\ldots,x_\mu)\in \mathbb{Z}[x_1,\ldots,x_\mu]$, $ \deg \Omega^{(\lambda)}\leq M$, and for every $p\in\mathbb{N}^*$, a polynomial $\Psi^{(\lambda)}_p(x_1,\ldots,x_{\mu+1})\in \mathbb{Z}[x_1,\ldots,x_{\mu+1}]$, $ \deg \Psi^{(\lambda)}_p\leq pM$, such that for every $y_0=\displaystyle\sum_{n\geq 1}c_nx^n$, $c_1\neq 0$, algebraic with vanishing polynomial of degrees bounded by $d_x$ in $x$ and $d_y$ in $y$, there exists $\lambda\in \Lambda$ such that for every $p\in\mathbb{N}^*$:
$$c_{\mu+1+p}=\displaystyle\frac{\Psi^{(\lambda)}_p(c_1,\ldots,c_{\mu+1})}{\Omega^{(\lambda)}(c_1,\ldots,c_\mu)^p}.$$
\end{coro}

\begin{demo}
Let $y_0=\displaystyle\sum_{n\geq 1}c_nx^n$, $c_1\neq 0$, be algebraic with vanishing polynomial of degree bounded by $d_x$ in $x$ and $d_y$ in $y$. According to the Theorem \ref{theo:wilc}, there is a finite set $\Lambda$ and for every $\lambda\in\Lambda$, polynomials $a^{(\lambda)}_{i,j}(x_1,\ldots,x_N)\in \mathbb{Z}[x_1,\ldots,x_N]$ such that:
\begin{equation}\label{equ:deg-coeff}
P^{(\lambda)}=\displaystyle\sum_{i\leq d_x,j\leq d_y}a_{i,j}^{(\lambda)}(c_1,\ldots,c_N)x^iy^j
\end{equation} 
is a vanishing polynomial for $y_0$ for a certain $\lambda\in\Lambda$. Enlarging the finite set $\Lambda$ by indices corresponding to the various $\displaystyle\frac{\partial^k P}{\partial y^k}$,$k=1,\ldots,d_y-1$, we can assume that there is $\lambda$ such that $y_0$ is a simple root of $P^{(\lambda)}$. So the coefficients of $y_0$ can be computed as in the Corollary \ref{coro:FS}. More precisely, for any $p\in\mathbb{N}^*$:
$$ c_{\mu+1+p}= \displaystyle\sum_{q=1}^p \displaystyle\sum_{S\in I_q} \displaystyle\sum_{T_S\in J_S}\displaystyle\frac{m_{S,T_S}}{\omega_0^p} \omega_0^{p-q}A^S C^{T_S}$$
where $I_q=\left\{(s_{i,j})\ |\ |S|=q,\ \|S\|_2\geq q-1\right\}$, $$J_S=\left\{(t_{S,i})\ |\ |T_S|=\|S\|_2-q+1,\ \|T_S\|=p+qi_\mu - (q-1)(\mu+1)-\|S\|_1\right\}$$ and $m_{S,T_S}\in \mathbb{Z}$. Note that $C=(c_1,\ldots,c_{\mu+1})$ and $A=(a_{i,j})$. It suffices to bound the degrees of the numerator and denominator in the terms of (\ref{equ:deg-coeff}). By the Theorem \ref{theo:wilc}, $\deg a_{i,j}^{(\lambda)}\leq M-d_y$. So by the Remark \ref{rem:omega_0}, we deduce that $\deg\omega_0\leq M$. The degree $d_{q,S}$ of a term $\omega_0^{p-q}A^S C^{T_S}$ is bounded by:
$$(p-q)(M-d_y)+|S|(M-d_y)+|T_S|=(p-q)(M-d_y)+q(M-d_y)+\|S\|_2-q+1.$$
But, $\|S\|_2\leq qd_y$ and $1\leq q\leq p$. So we get that:
$$d_{q,S}\leq p(M-d_y)+qd_y-q+1\leq pM$$
\end{demo}






\newcommand{\etalchar}[1]{$^{#1}$}
\def\cprime{$'$} \def\cprime{$'$}
\providecommand{\bysame}{\leavevmode\hbox to3em{\hrulefill}\thinspace}
\providecommand{\MR}{\relax\ifhmode\unskip\space\fi MR }
\providecommand{\MRhref}[2]{%
  \href{http://www.ams.org/mathscinet-getitem?mr=#1}{#2}
}
\providecommand{\href}[2]{#2}


%
\end{document}